\documentclass[12pt]{amsart}
\usepackage{amsmath,amssymb,amsfonts,amscd,graphicx}

\usepackage{chngcntr}
\usepackage{apptools}

\AtAppendix{\counterwithin{theorem}{section}}

\usepackage{array,multirow, hyperref}
\usepackage{color}

\usepackage{framed}

\usepackage{hyperref}
\hypersetup{
    colorlinks = true,
    linkcolor = {black},
   citecolor  = {black},
}

\usepackage{tikz}
\usetikzlibrary{arrows}
\usepackage{tikz-cd}
\usetikzlibrary{cd}

\DeclareMathOperator{\SO}{SO}
 \def\H{\mathbb{H}}

\usepackage{enumitem}
\usepackage{mathptmx}
\usepackage{amsthm}
\usepackage{amssymb}
\usepackage{tikz}
	\usetikzlibrary{decorations.pathreplacing}
	\usetikzlibrary{patterns}

\newcommand{\Q}{\mathbb{Q}}
\newcommand{\R}{\mathbb{R}}
\newcommand{\Z}{\mathbb{Z}}

\usepackage{xcolor}

\usepackage{changes}

\makeatletter
\newcommand{\addresseshere}{%
  \enddoc@text\let\enddoc@text\relax
}
\makeatother

\newtheorem{theorem}{Theorem}
\newtheorem{lemma}[theorem]{Lemma}

\newtheorem{proposition}[theorem]{Proposition}
\newtheorem{corollary}[theorem]{Corollary}
\theoremstyle{definition}

\theoremstyle{remark}
\newtheorem{remark}[theorem]{Remark}
\numberwithin{equation}{section}
\theoremstyle{plain}

\theoremstyle{definition}

\newcommand{\co}{\colon\thinspace}

\newcommand{\nl}{\hfil\break}

\usepackage{anyfontsize}

\theoremstyle{definition}

\newtheorem*{defin*}{Definition}

\theoremstyle{plain}

\usepackage{euscript}

\begin{document}


\title[]{\fontsize{12}{12}\selectfont  FILLING LINKS AND SPINES IN ${\mathbf 3}$-MANIFOLDS}

\author[]{\fontsize{12}{12}\selectfont Michael Freedman and Vyacheslav Krushkal\\
with Appendix by Christopher J. Leininger and Alan W. Reid}

\address{Michael Freedman, Microsoft Research, Station Q, and Department of Mathematics,
University of California, 
Santa Barbara, CA 93106, USA}
\email{mfreedman@math.ucsb.edu}

\address{Vyacheslav Krushkal, Department of Mathematics, University of Virginia,
Charlottesville, VA 22904, USA}
\email{krushkal@virginia.edu}

\begin{abstract} We introduce and study the notion of filling links in $3$-manifolds: a link $L$ is filling in $M$ if for any $1$-spine $G$ of $M$ which is disjoint from $L$, $\pi_1(G)$ injects into $\pi_1(M\smallsetminus L)$. A weaker ``$k$-filling'' version concerns injectivity modulo $k$-th term of the lower central series. For each $k\geq 2$ we construct a $k$-filling link in the $3$-torus. The proof relies on an extension of the Stallings theorem which may be of independent interest. We discuss notions related to ``filling'' links in $3$-manifolds, and formulate several open problems.
The appendix by 
C. Leininger and A. Reid establishes the existence of a
filling hyperbolic link
in any closed orientable $3$-manifold with $\pi_1(M)$ of rank $2$.
\end{abstract}

\begingroup
\def\uppercasenonmath#1{} % this disables uppercasing title
\let\MakeUppercase\relax % this disables uppercasing authors
\maketitle
\endgroup

\counterwithout{equation}{section}

\section{Introduction}

 There is an old theme in $3$-manifold topology that links or even knots can be built inside a general $3$-manifold which in one sense or another are as robust as an embedded  $1$-complex can be.\footnote{Links and $1$-complexes will be assumed to be PL embedded in a $3$-manifold.}  Here are several examples (in historical order): 
 
 \begin{enumerate}[leftmargin=.7cm,itemsep=3pt,parsep=2pt]
\item[1.]{Bing’s theorem \cite{Bing} that a closed $3$-manifold $M$ is diffeomorphic to $S^3$ if and only if every knot $K$ in $M$ is  contained (``engulfed'') in a $3$-ball.}
\item[2.]{``Disk busting curves'': Theorem \cite{M82}: for any compact $3$-manifold $M$ there is a knot $K$ in $M$ so that every essential sphere or disk meets $K$.  }
\item[3.]{Suppose $X\longrightarrow B\longrightarrow V$ is a manifold $(X)$-bundle over a $3$-manifold $V$. Then there is a link $L$ in $V$ so that $B$ restricted to $V\smallsetminus L$ admits a flat topological connection, i.e. a foliation of class $C^1$ transverse to the fibers \cite{M18, F20}. }
 \end{enumerate}
 
In all three examples the theorem becomes trivial if instead of considering a knot or link, we replace those words in the statement with an embedded $1$-complex, since the $1$-complex may be chosen to be a spine of a Heegaard handlebody.  So the pattern seems to be that by knotting and linking one can make $1$-submanifolds nearly as ``filling'' as a $1$-complex: hard to engulf, hard to avoid, and hard to flatten over.  Our first intent was to extend this theme to higher dimensions, but the natural questions have proved difficult and we just state a few of them for further thought:

\begin{enumerate}[leftmargin=.9cm,itemsep=3pt,parsep=2pt]
\item[Q1.]{ If $M$ is a smooth homotopy $4$-sphere with the property that every smoothly embedded surface (or perhaps just  every $2$-sphere?) lies in a ball, is M invertible?}

\item[Q2.]{ Is there an analog of (2) in high dimensions (being careful with respect to homotopy spheres)?}

\item[Q3.]{ In higher dimensions can such bundles be flattened over the complement of a codimension $2$ link, or even of a Cantor set?}
 \end{enumerate}

In this paper we consider only $3$-manifolds and attempt to add a $4$th example of ``filling'' by  a link or knot. To formulate the question, we start by discussing the relevant versions of the notion of a {\em spine} of a $3$-manifold $M$. We always consider $1$-complexes up to $I$-$H$ moves (also known as Whitehead moves) 
$ 
\vcenter{\hbox{\includegraphics[scale = 0.2]{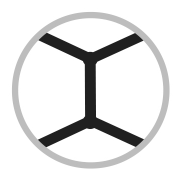}}} 
\leftrightarrow \vcenter{\hbox{\includegraphics[scale = 0.2]{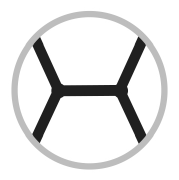}}}
$; considering their regular neighborhoods we also refer to them as handlebodies.

\begin{enumerate}[leftmargin=.7cm,itemsep=3pt,parsep=2pt] \label{spine def}
\item[(a)]{(The “rank” definition) A spine is a $1$-complex in $M$ of least first Betti number, surjecting onto $\pi_1(M)$.}
\item[(b)]{ A spine is a handlebody  in $M$ with the property that it is onto on $\pi_1$ but no smaller handlebody obtained from compression of a non-separating disk is onto.}
 \end{enumerate}

The following question may be formulated for each of the definitions above.

\begin{enumerate}[leftmargin=.9cm,itemsep=3pt,parsep=2pt]
\item[Q4.]{ 
 Given a compact $3$-manifold $M$, is there a link $L$ in $M$ so that whenever $G$ is a spine in $M$ and $G$ is disjoint from $L$, then $\pi_1 (G) \longrightarrow  \pi_1(M\smallsetminus L)$ is injective?
}
 \end{enumerate}

If there is such a link $L$, we call it {\em filling} in $M$. 
As we prove below any $1$-spine for $M$, if allowed to play the role of $L$ would have this injectivity property. So as in the first three examples the question is whether a link can ``do the work of'', or ``be as filling as'' a $1$-complex. This question and its relatives actually seem difficult to get hold of and we are only able establish some partial results. We present them in the hope that others will find these ``filling'' questions of interest.

We emphasize that that the term {\em spine} throughout the paper will refer either to a specific PL embedding of a $1$-complex into $M$, or to its regular neighborhood in $M$. One may consider spines up to equivalence: $I$-$H$ moves and isotopy in case (a), and either isotopy or homotopy of handlebodies in case (b). 
To check the filling condition using any of the definitions, it suffices to consider one representative from each equivalence class, with respect to the relevant equivalence relation ($I$-$H$ moves, isotopy, homotopy).
Considering spines up to isotopy or homotopy gives the same answer as to whether a link (or more generally a 1-complex) is or is not filling. The reason is that the difference between homotopy and isotopy is finger moves. And these finger moves can miss the link or 1-complex $L$ (that we are trying to decide if it is filling.)  These finger moves in the complement of $L$ do not change maps on $\pi_1$ so we get the same answers for isotopic and homotopic spines, when it comes to deciding if $L$ is filling.

We start with $\pi_1$-injectivity in the complement of a {\em $1$-complex}, the case that is much easier than the complement of a {\em link}, similarly to examples 1 -- 3 in the beginning of the introduction. In the following lemma, the $1$-complex that is shown to be filling is a spine of a Heegaard handlebody:
 
\begin{lemma} \label{spine lemma}\sl
Let $M=H\cup  H^*$ be a Heegaard decomposition. Then for each definition (a), (b), whenever $G$ is a spine in $M$ and $G$ is disjoint from $H^*$, the map
$\pi_1 (G)\rightarrowtail \pi_1(M\smallsetminus  H^*)$ induced by inclusion is injective.
\end{lemma}

{\em Proof.} Let $\pi$ be the image of $\pi_1(G)$ in $\pi_1(M\smallsetminus H^*)\cong \pi_1(H)$.  Being a subgroup of a free group, $\pi$ is free. Also ${\rm rank}({\pi})$ is less than or equal ${\rm rank}(\pi_1(G))$ since the map  $\pi_1(G) \longrightarrow \pi$ is onto.  If that map has a kernel, by the Hopfian property of free groups ${\rm rank} (\pi) < {\rm rank}(\pi_1(G))$, contradicting minimal rank in the context of definition (a).

Now consider a spine satisfying definition (b). 
As in the previous paragraph, $\pi_1(G) \longrightarrow \pi$ is onto, implying ${\rm rank} (\pi) < {\rm rank}(\pi_1(G))$. By \cite[Theorem 3.3]{MKS66}, there is a
basis $B$ for $\pi_1(G)$ so that a subset of $B$ generates ${\rm ker}[\pi_1(G)\longrightarrow \pi]$. 

\begin{proposition} \label{prop} \sl Let $G$ be a handlebody of genus $g>0$ and assume that $b \in \pi_1(G)$ is a free basis element. Then $G$ contains a properly embedded disk $D$ so that $$\pi_1(G)=\langle b\rangle *\pi_1(G\smallsetminus D).$$
\end{proposition}

{\em Proof of Proposition \ref{prop}}: Stated geometrically, the hypothesis says that there is a map $f\co G\longrightarrow W$ to a wedge of $g$ circles inducing an isomorphism on $\pi_1$ and taking $b$ to the generator of the first $S^1$ factor $S$ of $W$. Let $p \in S$ be a point on $S$ different from the basepoint, and $P \subset G$ be the transverse preimage, $f^{-1}(p)$, of $p$ in $G$. By standard $3$-manifold  techniques, we may assume that $P$ is an incompressible surface in $G$. (If there is an essential loop $\gamma$ on $P$ which is nullhomotopic in $G$, consider its null homotopy $J$. Choose an innermost essential circle $\alpha$ in $J$; $\alpha$ bounds a singular disk in $G\smallsetminus P$, which by Loop theorem/Dehn’s lemma may be replaced with an embedded disk $E \subset G\smallsetminus P$ with essential boundary in $P$. Compressing along $E$, and considering a homotopy of $f$ giving rise to this ambient surgery of the point preimage, reduces the complexity of $P$. This process must eventually make $P$ incompressible.) Since $f$ is an isomorphism on $\pi_1$ and $P$ is the preimage of a point it follows that  the map induced by inclusion: $\pi_1(P)\longrightarrow \pi_1(G)$ is the trivial map. So incompressibility implies that $P$ is a disjoint union of disks.  $P$ is Poincar\'{e} dual to the homology class $[b]$ so it follows that at least one of these disks $D$ is homologically essential, i.e. non-separating. This is the $D$ claimed in the proposition.
\qed

By Proposition \ref{prop},  a free generator in the kernel corresponds to a non separating reducing disk in $G$. Cutting along that disk gives a sub-handlebody still generating $\pi_1(M)$, contradicting minimality.
This concludes the proof of Lemma \ref{spine lemma}.
\qed

Now we turn to the main problem, the existence of filling links. While the notion of ``filling'' is interesting for both definitions (a) and (b), in this paper we will focus on the slightly simpler rank definition (a). Logically it is harder to show that a link is filling with respect to (b), since this definition allows more “spines” which need to all inject: a spine satisfying (b) might not have the least first Betti number among all $1$-complexes surjecting onto $\pi_1(M)$. As we see below, the analysis in the easier case (a) is already quite subtle.

One may find a {\em knot} $K$ giving  $\pi_1$-injectivity for a {\em fixed} embedding of a spine as follows. Consider a minimal genus Heegaard decomposition $M^3=H\cup  H^*$, and let $K$ be a diskbusting curve in the handlebody $H^*$. If $\pi_1 (H)\longrightarrow \pi_1(M\smallsetminus K)$ had kernel, by the loop theorem there would be a compressing disk in $ H^*$ disjoint from $K$, a contradiction.

The problem of finding a filling link $L$ (so an {\em arbitrary} embedding of a spine is $\pi_1$-injective in $M\smallsetminus L$) has a trivial solution for genus one $3$-manifolds.
In this case, $M^3=H\cup H^*$ where $H, H^*$ are solid tori, and a filling knot is given by the core circle of $H^*$.

The appendix by 
C. Leininger and A. Reid establishes the existence of a
filling hyperbolic link
in any closed orientable $3$-manifold with $\pi_1(M)$ of rank $2$. In particular, it follows that a closed orientable $3$-manifold with finite fundamental group contains a filling link, see Corollary \ref{finite}.

In the main body of this paper we focus on the case of the $3$-torus $T^3$, although it seems likely that the approach (which involves equivariant homological analysis in the universal cover of $M^3$) should work for hyperbolic $3$-manifolds as well.

The problem of analyzing $\pi_1$-injectivity of {\em any} embedding of a spine $G$ in the complement of a given link in the $3$-torus $T^3$ turned out to be quite subtle. The kernel of $\pi_1 (G)\longrightarrow \pi_1 (T^3)$ is the commutator subgroup of the free group $\pi_1 (G)$ on $3$ generators. A standard classical tool for showing injectivity of maps of the free group is the Stallings theorem, see Section \ref{background sec}.
However the Stallings theorem does not directly apply in our context. Specifically, it does not apply to the map $\pi_1 (G)\longrightarrow \pi_1(T^3\smallsetminus L)$ because it is not surjective on second homology, and it is injective, rather than an isomorphism, on $H_1$. The complexity of the problem reflects the fact that the image of the map on $\pi_1$ depends on the embedding of the handlebody $G$. One may attempt to apply the Stallings theorem to the map $[\pi_1 (G), \pi_1 (G)]\longrightarrow K$, where $K$ is the kernel $\pi_1(T^3\smallsetminus L)\longrightarrow \pi_1(T^3)$ for a suitable choice of $L$. But injectivity of the infinitely generated first homology of the commutator subgroup is hard to establish when the embedding $G\longrightarrow T^3\smallsetminus L$ changes by an arbitrary homotopy, as we discuss below. 
We are able to find links with a weaker property. We call a link $L\subset M$ ``$k$-filling'' if the injectivity in Q4 above holds modulo the $k$th term of the lower central series, $\pi_1 (G)/(\pi_1 (G))_k \rightarrowtail   \pi_1(M\smallsetminus L)/\pi_1(M\smallsetminus L)_k$.

\begin{theorem} \label{k-filling thm} \sl
For any $k\geq 2$ there exists a $k$-filling link in $T^3$.
\end{theorem}

To prove this theorem we give an extension of the Stallings theorem using powers of an augmentation ideal, which applies uniformly to all embeddings $G\longrightarrow T^3\smallsetminus L$, where the conclusion holds modulo a given term of the lower central series.
Powers of the augmentation ideal of a group ring have been classically studied in low-dimensional topology, cf. \cite{CG83, S75}, and are related to the lower central series of the group. We analyze a different aspect of the theory, focusing on the connection between the lower central series of a group $\pi$ and powers of the augmentation ideal of the group ring of $H_1(\pi)$.

It is interesting to note the similarity of Theorem \ref{k-filling thm} with the current state of knowledge of the topological $4$-dimensional surgery theorem for free non-abelian groups (see \cite{FQ90}). In the underlying technical statement, the disk embedding conjecture, one considers disks up to homotopy fixing the boundary, and the question is whether the map on ${\pi}_1$ can be made {\em trivial}. One can solve the problem modulo any term of the lower central series, but the question itself is open. (See \cite{FK20} for recent developments.)

In an October 2021 arXiv posting \cite{Stagner}, 
William Stagner proved that $3$-manifolds of rank $3$ have filling links. In particular, his work gives a genuine filling link in the $3$-torus, thus strengthening our Theorem \ref{k-filling thm}. The result of \cite{Stagner} relies on methods of hyperbolic geometry, very different from our approach.

The proof of Theorem \ref{k-filling thm} follows from Lemmas \ref{augmentation ideal lemma}, \ref{links lemma}, \ref{augm lcs} given in the following sections. Section \ref{FM and powers} and Lemma \ref{augmentation ideal lemma} set up the equivariant homological framework for analyzing the effect of homotopies of a spine in terms of powers of the augmentation ideal. Lemma \ref{links lemma} and its analogue for $T^3$ in section \ref{3-torus section} construct links satisfying the conditions of Lemma \ref{augmentation ideal lemma}. Lemma \ref{augm lcs} in section \ref{sec: exact seq} gives an extension of the Stallings theorem, relating powers of the augmentation ideal to the lower central series, needed to complete the proof of the theorem.

\section{The relative case: $T^2\times I$} \label{FM and powers}

\subsection{Notation and background} \label{background sec}
The coefficients of homology groups are set to be $\mathbb Z$ throughout the paper. We start by recalling the Stallings theorem. Given a group $A$, its lower central series is defined inductively by $A_1=A$, $A_k=[A_{k-1}, A]$; $A_{\omega}=\cap_{k=1}^{\infty} A_k$.

\begin{theorem}[Stallings' theorem \cite{S65}] \label{Stallings theorem}  \sl Let $f\co A\longrightarrow B$ be a group homomorphism inducing an isomorphism on $H_1$ and an epimorphism on $H_2$. Then $f$ induces an isomorphism $A/A_k \longrightarrow B/B_k$ for all finite $k$, and an injective map $A/A_{\omega}\longrightarrow B/B_{\omega}$.
\end{theorem}

Dwyer \cite{D75} extended the theorem, relaxing the surjectivity to be onto $H_2$ modulo the $k$-th term of the {\em Dwyer filtration}:
$$\phi_n(A)\; =\; {\rm ker}[H_2(A)\longrightarrow H_2(A/A_n)]
$$

Assuming that $f\co A\longrightarrow B$ is an isomorphism on $H_1$, the result of \cite{D75} is that $f$ induces an isomorphism $A/A_{k+1} \longrightarrow B/B_{k+1}$ if and only if it induces an epimorphism
$H_2(A)/\phi_k(A)\, \longrightarrow \, H_2(B)/\phi_k(B)$.

\subsection{Relative spine and equivariant homology in the universal cover} 
\begin{figure}[ht]
\includegraphics[height=3.5cm]{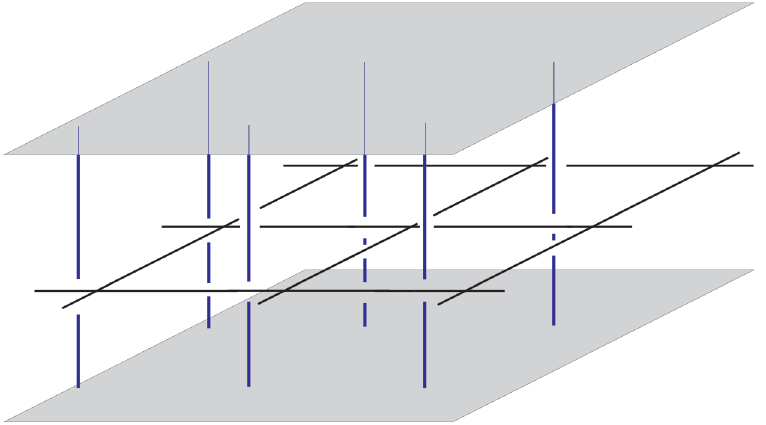}
\caption{The preimage $\widetilde G$ in the universal cover  $\R^2\times I$ of the standard relative spine $G=(\{*\}\times I)\cup(T^2\times \partial I)$ consists of the top and bottom shaded panels union the vertical line segments. The mid-level horizontal grid is the preimage of the dual spine $S^1\vee S^1\subset T^2\times \{1/2\}$.}
\label{fig:Relative spine}
\end{figure}

In this section we consider the relative case where the construction of $k$-filling links is easier to describe, $M=T^2\times I$. Fix the standard ``relative spine'' 
$G=(\{ * \}\times I)\cup(T^2\times \partial I)$, 
and the dual spine $G^*=S^1\vee S^1\subset T^2\times \{1/2\}$. Their preimages in the universal cover are illustrated in figure \ref{fig:Relative spine}. 

An analogue of lemma \ref{spine lemma} shows that for any proper embedding of the line segment $\{ *\}\times I$ into $((T^2\times I)\smallsetminus G^*, T^2\times \partial I)$ the fundamental group $\Z^2 *\Z^2$ of the resulting relative spine injects into $\pi_1(T^2\times I\smallsetminus G^*)$. This may be seen using the Stallings theorem. Note that $\pi_1 (G)\cong \Z^2 *\Z^2$ is a right-angled Artin group, and it is also a free product of surface groups; these classes of groups are residually nilpotent \cite{DK92}, \cite{CO98}. (Recall that the Stallings theorem, under the assumptions stated in Theorem \ref{Stallings theorem}, gives injectivity modulo the $\omega$-term of the lower central series. As noted above, $\pi_1 (G)$ is residually nilpotent, thus $(\pi_1 (G))_{\omega}$ is trivial, and therefore the Stallings theorem implies injectivity of $\pi_1 (G)$.)

\begin{figure}[ht]
\includegraphics[height=3.5cm]{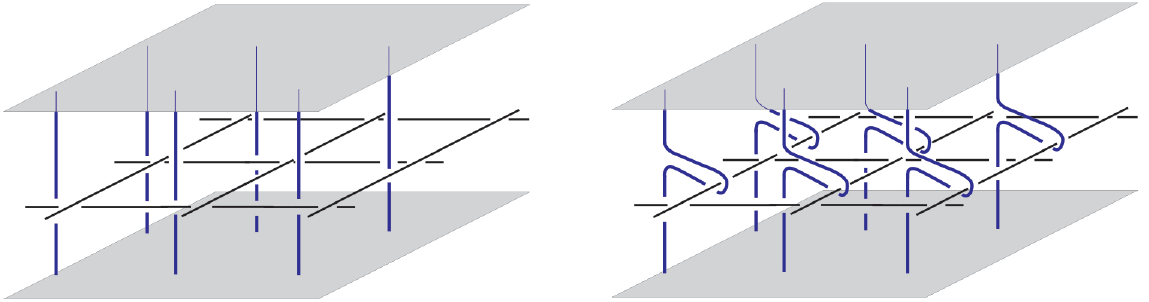}
\caption{An example of a link $L\subset T^2\times (0,1)$, and a finger move.}
\label{fig:finger move}
\end{figure}

The goal is to analyze the map on $\pi_1$ induced by inclusion when $G^*$ is replaced by a link. This is illustrated in figure \ref{fig:finger move} where the link $L\subset T\times (0,1)$ is obtained by ``resolving'' the dual spine $G^*$ into two disjoint essential circles. Different embeddings of the vertical interval are related by homotopies that may pass through link components and may be thought of as finger moves; the equivariant lift of  an elementary finger move is illustrated on the right in figure \ref{fig:finger move}.

We will continue denoting by $G$ the standard embedding of the relative spine into $T^2\times I$, and by $\widetilde G$ its preimage in the universal cover. The notation $G', \widetilde G'$ will be used when the embedding of the vertical arc is arbitrary, i.e. related to the standard embedding by finger moves.
The fundamental group of $\widetilde G$, and also of $\widetilde G'$ is $$K:={\rm ker}[\Z^2 *\Z^2\longrightarrow \Z^2].$$

Let $L$ be a link in $T^2\times (0,1)$ whose components are all essential in $\pi_1(T^2)$, and let $\widetilde L$ denote its preimage: a $\Z^2$-equivariant collection of lines in the universal cover. In our examples the links will consist of geodesics in $T^2$ shifted to different levels in the interval $(0,1)$. In this case the preimage $\widetilde L$ consists of a disjoint collection of straight lines. Considering a radial Morse function, the complement $\R^2\times I\smallsetminus \widetilde L$ is seen to have a handle decomposition with infinitely many $1$-handles (one for each line) and no $2$-handles, so  $\pi_1(\R^2\times I\smallsetminus \widetilde L)$ is free. More general links would require an extension of the group theoretic analysis developed in this paper.
The starting point is to analyze the injectivity of the map $\alpha$ in the commutative triangle
$$
\centering
\begin{tikzcd}
 \pi_1( G) \arrow[r, "\alpha"] \arrow[dr, "\beta"']  & \pi_1( T^2\times I\smallsetminus L) \arrow[d,"\gamma"'] \\
& \pi_1(T^2\times I)  & 
\end{tikzcd}
$$

where all maps are induced by inclusions. Here ${\alpha}$ is induced by the {\em standard} embedding $G\subset T^2\times I\smallsetminus L$; the general case of $G'$ is discussed in section \ref{fm subsec}. It is clear that ${\rm ker}[\alpha]$ equals the kernel of the map ${\rm ker}[\beta]\longrightarrow {\rm ker}[\gamma]$. Here ${\rm ker}[\beta]=K$, and ${\rm ker}[\gamma]\cong \pi_1(\R^2\times I\smallsetminus \widetilde L)$. Therefore the focus is on the map
\begin{equation} \label{incl map eq} K\longrightarrow \pi_1(\R^2\times I\smallsetminus \widetilde L).
\end{equation}

Denote by $J$ the first homology of $\widetilde G$, $J=K/[K,K]$, and let $H$ denote $H_1(\R^2\times I\smallsetminus \widetilde L)$. Since $\pi_1(\R^2\times I\smallsetminus \widetilde L)$ is a free group, if $J\longrightarrow H$ were injective, the Stallings theorem would imply that the map (\ref{incl map eq}) is injective (for the standard spine $G$). We will give a construction of links $L$ satisfying weaker injectivity, modulo powers of the augmentation ideal (see below), for any embedding $G'\subset T^2\times I\smallsetminus L$.

Denote the map $J\longrightarrow H$ by $Lk$. The group $H$ is generated by meridians $m(l)$, one for each line $l$ in $\widetilde L$. The map $Lk$ is given by $\Z^2$-equivariant linking, sending a $1$-cycle $c$ in $\widetilde G$ to a linear combination of meridians $\sum_i\,  a_i\,  m(l_i)$, where the coefficient $a_i\in\Z[\Z^2]$ is the linking ``number'' of $c$ and $l_i$. Since there is a single generator $m(l)$ for each line $l$, when there is no risk of confusion we will write
\begin{equation} \label{m l} 
Lk(c)\, =\, \sum_i \, a_i \,  l_i.
\end{equation}
An example calculating the linking map is given below.

\begin{figure}[ht]
\includegraphics[height=3.5cm]{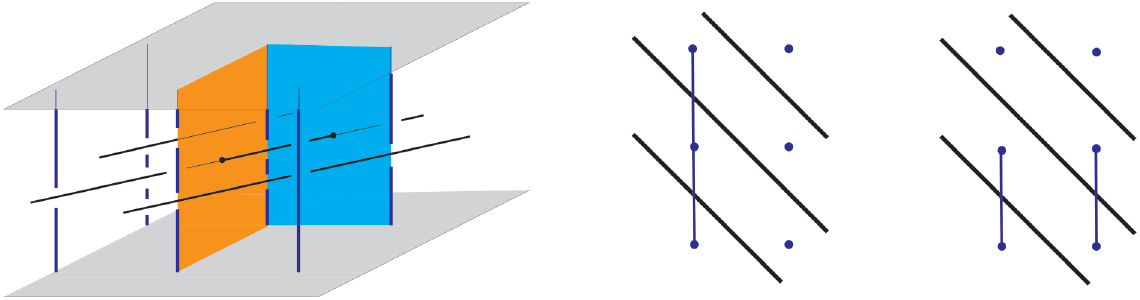}
\scriptsize
\put(-316, 20){$P_x$}
\put(-265, 30){$P_y$}
\put(-395,0){(a)}
\put(-161,23){$P_x$}
\put(-169,57){$-yP_x$}
\put(-58,23){$P_x$}
\put(-34.5,23){$-xP_x$}
\put(-180,0){(b)}
\put(-70,0){(c)}
\put(-116,2){$l$}
\put(-14,2){$l$}
\caption{\nl (a): Plaquettes $P_x, P_y$ generating $J$ over $\Z[\Z^2]$. \nl (b), (c): Projection onto $\R^2$; dots represent the preimage of the edge $\{ * \}\times I$ of the relative spine.}
\label{fig:relative link}
\end{figure}

As a module over $\Z[\Z^2]$, $J$ is generated by the boundaries of two vertical ``plaquettes'', denoted $P_x$ and $P_y$ in figure \ref{fig:relative link}a. We will think of elements of $J$ as linear combinations of these plaquettes, with coefficients in $\Z[\Z^2]$. The translations in the directions perpendicular to $P_x, P_y$ are denoted respectively by $x,y$. Note that the relation 
\begin{equation} \label{relation eq}
(1-x) P_x+(1-y)P_y\, =\, 0
\end{equation}

holds in J. 

Figure \ref{fig:relative link} illustrates the case when the link $L$ has a single component, the $(1, 1)$-curve in the torus $T^2\times \{ 1/2\}$. In this case the two translations act the same way on $\widetilde L$: for any line $l$,  $xl=yl$. Figures \ref{fig:relative link} (b, c) show the projection onto $\R^2$ of two elements of $J$: $(1-y)P_x$, $(1-x)P_x$. Denoting the line intersecting the plaquette with boundary $ P_x$ by $l_0$, we have
\begin{equation}\label{meridian line}  Lk((1-y)P_x)\, =\, (1-y)\, l_0, \; \; Lk((1-x)P_x)\, =\, (1-x)\,  l_0.
\end{equation}

(In the equation above, following \eqref{m l} and the sentence preceding it, for brevity of notation we use $l_0$ also to denote the meridian to the line.)  
Since $(1-y)\, l_0 =(1-x)\, l_0$ in this example, the map $J\longrightarrow H$ is certainly not injective. 
Since the map $Lk$ turned out to have kernel in this example, it is no surprise to observe that the boundary curve of $ P_x\cup  P_y$ in Figure \ref{fig:relative link} (a) is null homotopic in $\R^2\times I\smallsetminus \widetilde L$.

\subsection{Finger moves and the kernel of the linking map} \label{fm subsec}
Next we extend this discussion to the case of an arbitrary embedding $G'\subset T^2\times I\smallsetminus L$. 
Consider the cellular chain complex of the preimage $\widetilde G$ of the standard spine $G$.
Up to a homotopy equivalence where $\R^2\times \{ 0\}$ and $\R^2\times \{ 1\}$ are contracted to points, we will consider $\widetilde G$ as a $1$-complex.   
Since there are no $2$-cells in $\widetilde G$, the first homology $J$ will be identified with cellular $1$-cycles.
Let 
$j$ denote the inclusion of cellular $1$-cycles into cellular $1$-chains of $\widetilde G$, $j\co J\hookrightarrow C_1$. Here 
$C_1$ is generated by a single vertical line segment as a module over $\Z[\Z^2]$. 

The linking map $Lk'\co J'\longrightarrow H$ for any $G'$ is given by $Lk+F\circ j$, where $F\co C_1\longrightarrow H$ is the ``finger move'' map which measures the difference in $H$ between the standard embedding $\{ * \}\times I$ and its homotopic image in $G'$. It follows that a $1$-cycle $c$ is in the kernel of $Lk'$ if and only if
\begin{equation} \label{c in ker}
Lk(c)=-F(j(c))
\end{equation}

It is convenient to represent this using diagram (\ref{FM cd}) of $\Z^2$-equivariant maps. (Note that this diagram does {\em not} commute; the two maps are equal precisely on the kernel of $Lk'$.)

\begin{equation} \label{FM cd}
\centering
\begin{tikzcd}
  J \arrow[r,"Lk"] \arrow[d,"j"'] & H \\
 C_1 \arrow[ur,"-F"'] & 
\end{tikzcd}
\end{equation}

\subsection{Powers of the augmentation ideal}
Consider elements of ${\mathbb Z}[{\mathbb Z}^2]$ as Laurent polynomials in two commuting variables $x, y$. Let $I$ denote the augmentation ideal of ${\mathbb Z}[{\mathbb Z}^2]$. 
The following lemma provides a convenient tool for analyzing the injectivity of the linking map for an arbitrary spine, modulo powers of the augmentation ideal.

\begin{lemma} \label{augmentation ideal lemma}{\sl
If $$i_k\co I^k J/I^{k+1} J\longrightarrow  I^k H/I^{k+1}  H$$ is injective for some $k$, then for any relative spine $G'\subset T^3\smallsetminus L$, $$i'_k \co I^k J' /I^{k+1} J'\; \longrightarrow \; I^k H/I^{k+1}  H$$ is injective. Here $i_k, i'_k$ are the maps induced by the inclusions of $G, G'$ into $ T^3\smallsetminus L$.
}
\end{lemma}

{\em Proof}. 
The maps in (\ref{FM cd}) are equivariant over $\Z[\Z^2]$; consider
\begin{equation} \label{FM cdk}
\centering
\begin{tikzcd}
  I^k J/I^{k+1} J \arrow[r,"Lk_k"] \arrow[d,"j_k"'] & I^k H/I^{k+1} H \\
 I^k C_1/I^{k+1} C_1 \arrow[ur,"-F_k"'] & 
\end{tikzcd}
\end{equation}

The module $J$ is generated by the plaquette boundaries $P_x, P_y$, and their images under $j$ in $C_1$ are 
$$j(P_x)=(1-y) Z, \; j(P_y)=(1-x) Z, $$

where $Z$ denotes a generator (vertical line segment) of the
module $C_1$. 
Suppose $c\in I^k J$. Then 
$c=aP_x+bP_y$ for some $a,b\in I^k$, and
$$
j_k(c)\, =\, (a(1-y)+b(1-x)) Z\, \in \, I^{k+1}C_1.
$$

It follows that $j_k(I^k J)\subset I^{k+1} C_1$, so $j_k=0$, and therefore
$${\rm ker}[Lk_k+F_k\circ j_k]\; = \; {\rm ker}[Lk_k]. 
$$
It follows that the kernel is independent of the finger move map $F$; thus it is independent of the choice of $G'$.
\qed

\subsection{Construction of the links} \label{construction sec}
Given $k\geq 0$, in this section we construct links in $T^2\times I$ used in the proof of Theorem \ref{k-filling thm}.

\begin{lemma}\label{links lemma} \sl For any $k$ there exists a link $L_k\subset T^2\times I$ such that $i_j\co I^j J/I^{j+1} J\longrightarrow  I^j H/I^{j+1}  H$ is injective for all $1 \leq j\leq k$.
\end{lemma}

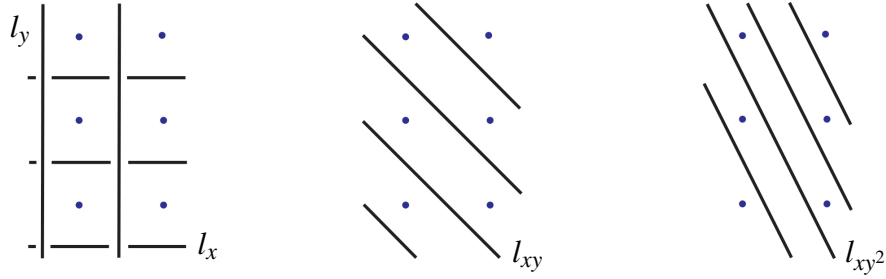
\begin{figure}[ht]
\includegraphics[height=3.5cm]{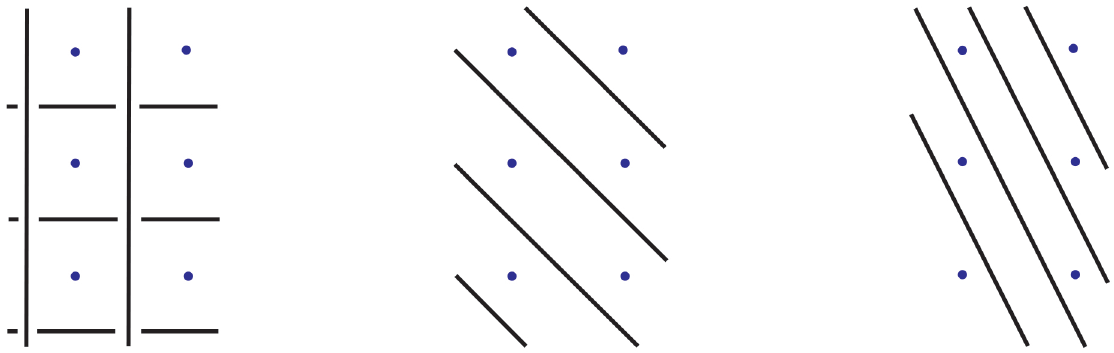}
\put(-249,3){$l_x$}
\put(-320,85){$l_y$}
\put(-130,0){$l_{xy}$}
\put(-4,0){$l_{xy^2}$}
\caption{The preimage of the curves in $T^2\times I$, used in the proof of lemma \ref {links lemma}, in the universal cover: projection of $\R^2\times I$ onto $\R^2$ is shown; as in figure \ref{fig:relative link} dots represent the preimage of the edge $\{ * \}\times I$ of the relative spine.}
\label{fig:Link1}
\end{figure}

{\em Proof.} To illustrate the idea of the proof, consider small values of $k$. For $k=0$, let $L_0$ be the $2$-component link consisting of a symplectic basis $(1,0), (0,1)$ curves on the torus, shifted to disjoint levels in $T^2\times I$. Two lines in their preimage in $\R^2\times I$ are denoted $l_x, l_y$, figure \ref{fig:Link1}.
In this case $J/IJ$ is $2$-dimensional, generated by the plaquettes $P_x, P_y$, figure \ref{fig:relative link}. The linking map $i_0\co J/IJ\longrightarrow H/IH$ is represented by the identity $2\times 2$ matrix over the integers. (Here the linking numbers are elements of $\Z[\Z^2]/I\cong \Z$.)

For $k=1$, $L_1$ is defined to be $L_0$ union the $(1,-1)$ curve on the torus, shifted to a different level in $T^2\times I$; a line in its preimage is denoted $l_{xy}$. (In our convention the slopes are negative, as shown in Figures \ref{fig:relative link}, \ref{fig:Link1}.) In this case $IJ/I^2J$ is $3$-dimensional, spanned by $(1-x) P_x, (1-y) P_x, (1-x)P_y$
(recall the relation (\ref{relation eq}) $(1-x) P_x+(1-y)P_y = 0$). 

$IH/I^2H$ is generated by $(1-x), (1-y)$ times the meridians to the lines $l_x, l_y, l_{xy}$. As in equation \eqref{m l}, we continue abbreviating to $l$ the notation $m(l)$ for the meridian. Note that $$(1-x) l_x=0, \, (1-y) l_y=0, \, (1-x)l_{xy}=(1-y)l_{xy}.
$$
Therefore $IH/I^2H$ is $3$-dimensional. 
The $3\times 3$ matrix 
representing the linking map $i_1\co IJ/I^2J\longrightarrow IH/I^2H$ is
$$
\begin{array}{c|ccc}
\multicolumn{4}{r}{ \hspace{1.9cm} {{{(1-x)l_y}}} \hspace{.5cm} {{{(1-y)l_x}}}     \hspace{.5cm} {{{(1-x)l_{xy}}}}\;\;\, \;}\\
\hline
  {{{(1-x)P_y}}} &  \hspace{.6cm} {1} & \hspace{1.4cm} {0}  & \hspace{.4cm} {1} \\
   {{{(1-y)P_x}}} & \hspace{.6cm} {0} & \hspace{1.4cm} {1} & \hspace{.4cm} {1} \\
    {{{(1-x)P_x}}} &  \hspace{.6cm} {0} & \hspace{1.4cm} {0} & \hspace{.4cm} {1}  \\
%\multicolumn{5}{c}{}
\end{array}
$$

The map $i_1$ is seen to be injective. 
The final example we consider is $k=2$. Using relation (\ref{relation eq}), consider  the generators $$(1-x)^2 P_x, (1-x)(1-y) P_x, (1-y)^2 P_x, (1-x)^2 P_y$$
of $I^2J/I^3 J$. Let $L_2$ denote the $4$-component link, obtained from $L_1$ by adding the $(1,-2)$ curve on the torus; a line in its preimage is denoted $l_{xy^2}$ in figure \ref{fig:Link1}. Like in the previous case, linking with $l_{xy}$ is defined modulo $x=y$.

The actions of $x, y$ on the new line are related by 
$xl_{xy^2}=y^2l_{xy^2}$.
Note the equality 
\begin{equation} \label{equality mod I2}
(1-y^2)=(1-y)(1+y)= 2(1-y) \;\;  {\rm in} \; \; I/I^2
\end{equation}
Using (\ref{equality mod I2}), for example the linking number of $(1-x)^2P_x$ and $l_{xy^2}$ in $I^2/I^3$ equals $(1-x)^2=(1-y^2)^2=4(1-y)^2;$ so the linking number of $(1-x)^2P_x$ and $(1-y)^2l_{xy^2}$ in the table below is $4$.

$$
\begin{array}{c|cccc}
\multicolumn{5}{r}{ \hspace{3cm} {{{(1-x)^2 l_y}}} \hspace{.5cm} {{{(1-y)^2 l_x}}}     \hspace{.5cm} {{{(1-y)^2 l_{xy}}}} \hspace{.5cm} {{{(1-y)^2 l_{xy^2}}}}\;\;\;\;\; \,\;}\\
\hline
{{{(1-x)^2 P_y}}} & \hspace{.6cm} {1} & \hspace{1.5cm} {0}  & \hspace{1.5cm} {1}  & \hspace{.4cm} {4}\\
  {{{(1-y)^2 P_x}}} & \hspace{.6cm} {0}    & \hspace{1.5cm} {1} & \hspace{1.5cm} {1}  & \hspace{.4cm} {1}\\
  {{{(1-x)(1-y) P_x}}} & \hspace{.6cm} {0} & \hspace{1.5cm} {0} & \hspace{1.5cm} {1} & \hspace{.4cm} {2} \\ 
   {{{(1-x)^2 P_x}}} & \hspace{.6cm} {0} & \hspace{1.5cm} {0} & \hspace{1.5cm} {1} & \hspace{.4cm} {{4}} \\
%\multicolumn{5}{c}{}
\end{array}
$$

The lower right $2\times2$ block of the matrix is $\scriptsize{\left(\!\!\begin{array}{cc}
   1  & 2 \\
    1 & 4
\end{array}\!\!\right)}$, a non-singular matrix. Thus $i_2$ is injective.

Continuing the inductive construction, the link $L_k\subset T^2\times I$ is defined to be $L_{k-1}$ union the $(1,-k)$ curve on the torus, shifted in the $I$ coordinate to be disjoint from $L_{k-1}$.
The group $I^kJ/I^{k+1}J$ is generated by $k+2$ elements: 
$$ \{ \, (1-x)^a(1-y)^b P_x, \; a+b=k,\; a,b\geq 0\, \} \; \; {\rm and}\; \; (1-x)^k P_y.
$$

Using the equality $1-y^j=j(1-y)$ in $I/I^2$, the linking number in $I^k/I^{k+1}$ of $(1-x)^a (1-y)^{k-a} P_x$,  and $l_{xy^j}$ equals \begin{equation} \label{calc}
(1-x)^a (1-y)^{k-a}=(1-y^j)^a (1-y)^{k-a}=j^a(1-y)^a (1-y)^{k-a}=j^a(1-y)^k.
\end{equation}
Therefore the linking matrix has an upper triangular block decomposition with the diagonal $2\times 2$ block $\scriptsize{\left(\!\!\begin{array}{cc}
   1  & 0 \\
    0 & 1
\end{array}\!\!\right)}$ and the $k\times k$
Vandermonde matrix $V$ whose $(m,n)$-th entry is $V_{m,n}=m^{n-1}$. Its determinant is non-zero, showing that $i_k\co I^k J/I^{k+1} J\longrightarrow  I^k H/I^{k+1}  H$ is injective.
This concludes the proof of lemma \ref{links lemma}.
\qed

\section{the $3$-torus} \label{3-torus section}
This section extends the construction from the relative case $T^2\times I$ above to $T^3$. 
After setting up the notation and giving details of the construction, we will outline the analogues of 
lemmas \ref{augmentation ideal lemma}, \ref{links lemma}.

\begin{figure}[ht]
\includegraphics[height=4.5cm]{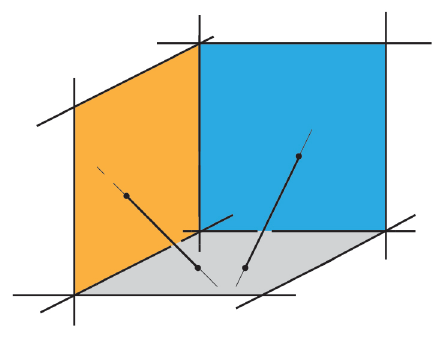}
{\scriptsize
\put(-105,94){$P_x$}
\put(-36,95){$P_y$}
\put(-54,32){$P_z$}
\put(-110,47){$l_{xz}$}
\put(-73,56){$l_{yz}$}
}
\caption{Lines and plaquettes.}
\label{fig:lines}
\end{figure}
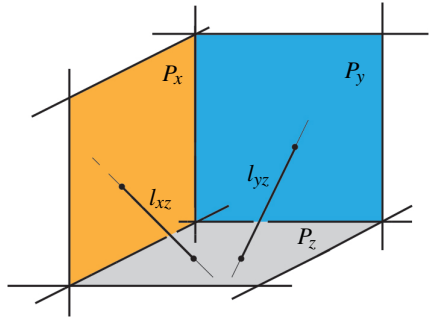

Let $T^3=H\cup  H^*$ be the genus $3$ Heegaard decomposition, and let $G, G^*$ denote the spines of the two handlebodies. Denote $F:=\pi_1 (G)$, the free group on three generators. The fundamental group of the preimage $\widetilde G$ in $\R^3$ is the commutator subgroup $F^{}_2$. Denote by $J$ the first homology of the cubical lattice $\widetilde G$ (the ``jungle gym''),  $J=F^{}_2/[F^{}_2,F^{}_2]$. As in the previous section, we think of generators of $J$ as plaquettes. As a module over $\Z[\Z^3]$, $J$ is generated by the boundaries of the plaquettes $P_x, P_y, P_z$, figure \ref{fig:lines}. Denoting the three covering translations by $x,y,z$, observe the relation
\begin{equation} \label{relation eq2}
(1-x) P_x+(1-y)P_y+(1-z)P_z\, =\, 0,
\end{equation}

in $J$, analogous to (\ref{relation eq}).

Consider a link $L\subset T^3 \smallsetminus G$ whose components are homologically essential in $T^3$. Its preimage, a $\Z^3$-equivariant collection of lines in ${\mathbb R}^3$, will be denoted $\widetilde L\subset {\mathbb R}^3$. Consider
$$K:= {\rm ker}\, [\, \pi_1(T^3\smallsetminus L)\longrightarrow \pi_1(T^3)\, ]. $$

$K\cong\pi_1({\mathbb R}^3\smallsetminus \widetilde L)$ is a free group. Denote by $H$ the first homology of ${\mathbb R}^3\smallsetminus \widetilde L$, $H=K/[K,K].$
Consider the map $Lk\co J\longrightarrow H$ induced by the inclusion $G\hookrightarrow T^3\smallsetminus L$. It is the $\Z[\Z^3]$-equivariant linking map, or in other words it is given by equivariant intersection between the plaquettes and the lines $\widetilde L$.  
Elements of $R:=\Z[\Z^3]$ are represented as Laurent polynomials in commuting variables $x, y, z$. Denote the augmentation ideal of $R$ by $I$. 

Let $i'\co G'\longrightarrow T^3\smallsetminus L$ be a spine homotopic to $G$. Such a homotopy may be thought of as finger moves of the three edges of $G$, intersecting the components of $L$ during the homotopy. Lifting this to the universal cover, we have a $\Z[\Z^3]$-equivariant map $F\co C_1(\widetilde G)\longrightarrow H$.
Here the chain group $C_1$ of $\widetilde G$ is a free module over $R$ of rank three. With this notation in place, the statement and the proof of lemma \ref{augmentation ideal lemma} 
hold without any changes.

Next we adapt the construction of links $L_k$ in lemma \ref{links lemma} to the setting of $T^3$. For $k=0$ it suffices to consider a $3$-component link $L_0$ obtained as a ``resolution'' of the $1$-spine $ G^*$ of $T^3$. Its preimage in $\R^3$, the lines $l_x, l_y, l_z$, pair $\delta_{i,j}$ with the plaquette generators of $J$, and the map $J/IJ\longrightarrow H/IH$ is an isomorphism for $L_0$.

Consider generators of $I^k/I^{k+1}$ of the form \begin{equation} \label{basis elements}
C^k_{a,b,c}:=(1-x)^a(1-y)^b(1-z)^c, \;\, a+b+c=k, \; a,b,c\geq 0. 
\end{equation}
It follows from equation (\ref{relation eq2}) that 
\begin{equation}
\{ C^k_{a,b,c}\, P_x,\;   C^k_{a,b,c}\, P_y, \;  C^k_{a,b,0}\, P_z    \}  
\end{equation}
is a basis for $I^kJ/I^{k+1}J$. Consider the $(2k+3)$-component link $L_k$ given by $L_0$ union $(1, 0, -j)$ and $(0, 1, -j)$ curves in $T^3$,  $j=1,\ldots, k\}$.
Their preimages are denoted
$l_{xz^j},\; l_{yz^j}$. The claim is that as in lemma \ref{links lemma},  $i_j\co I^j J/I^{j+1} J\longrightarrow  I^j H/I^{j+1}  H$ is injective for all $1 \leq j\leq k$.
The proof is an extension of the analysis in section \ref{construction sec}, which we outline next.

The linking pairing over $I^k/I^{k+1}$ between the basis elements $B_0$ of $I^kJ/I^{k+1}J$ of the form $\{ C^k_{0,b,c}\, P_x,\;   C^k_{a,0,c}\, P_y, \;  C^k_{a,b,0}\, P_z    \}$ and the lines $l_x, l_y, l_z$ is non-degenerate. (For example, the linking number of  $C^k_{0,b,c}\, P_x$ with $l_x$ equals $(1-y)^b(1-z)^c$ and with $l_y, l_z$ it is zero.) Since the linking pairing of $l_x, l_y, l_z$ is trivial with all other basis elements (denote them $B_{\neq 0}$) of $I^kJ/I^{k+1}J$, it suffices to analyze the linking pairing between $B_{\neq 0}$ and the lines $\{ l_{xz^j},\; l_{yz^j}\}$. 
It decomposes as the direct sum of the linking matrix for 
\begin{equation} \label{x eq}
\{ C^k_{a\neq 0,b,c}\, P_x, \, a+b+c=k,\; a,b,c\geq 0\} \; \; {\rm  and} \;\; \{ l_{xz^j}, j=1,\ldots, k\}
\end{equation}
and the analogous pairing for $\{ C^k_{a,b\neq 0,c}\, P_y\} $ and $\{ l_{yz^j}\}$. Consider (\ref{x eq}); the analysis for $y$ in place of $x$ is directly analogous.

The linking number (element of $I^k/I^{k+1}$) of $C^k_{a,b,c}\, P_x$ with $l_{xz^j}$ equals $(1-x)^a(1-y)^b(1-z)^c= j^a(1-y)^b(1-z)^{a+c}$, see (\ref{calc}). For different values of $b$, the linking numbers correspond to different basis elements of $I^k/I^{k+1}$. For a given value of $b$, the linking matrix of $\{C^k_{a,b,c}\, P_x, \, a=1,\ldots, k-c-b\}$ and  $\{(1-y)^b(1-z)^{a+c} l_{xz^j}, j=1,\ldots, k-c-b\}$ is the Vandermond matrix, as in the proof of lemma \ref{links lemma}, showing that $i_k$ is injective as claimed.

{\em Remark.} If a link $L$ is replaced with the dual spine $ G^*$ (as in the setting of lemma \ref{spine lemma}), the analysis above substantially simplifies. Recall the following result, an instance of a more general statement in terms of the torsion-free derived series from \cite{CH08}:

\begin{theorem} \label{CH thm}\cite[Corollary 2.2]{CH08} {\sl Suppose $F$ is a free group, $B$ is a finitely-related group, $\phi\co F \longrightarrow B$
induces a monomorphism on $H_1(-; \Q )$, and $H_2(B; \Q )$ is spanned by $B^{(n)}$-surfaces.
Then $\phi$ induces a monomorphism $F/F^{(n+1)} \subset B/B^{(n+1)}$.}
\end{theorem}

Here for a group $B$, its derived series is defined by $B^{(0)}=B, \, B^{(n+1)}=[B^{(n)},B^{(n)}]$. The notion of $B^{(n)}$-surfaces (maps of surfaces into $K(B,1)$ where the image on $\pi_1$ is in $B^{(n)}$ \cite[Definition 1.5]{CH08}) gives an analogue of the Dwyer filtration in the derived setting.

Any embedding $i\co G\longrightarrow T^3\smallsetminus  G^*$, induces an isomorphism on $H_1$, because both $\phi, \psi$ in the diagram below are isomorphisms:

\begin{equation} \label{spine eq}
\centering
\begin{tikzcd}
 H_1(G) \arrow[r, "\phi"] \arrow[d,"i_*"'] & H_1( T^3) \\
 H_1(T^3\smallsetminus  G^*) \arrow[ur, "\psi"'] & 
\end{tikzcd}
\end{equation}

Theorem \ref{CH thm} applies to $F=\pi_1(G), B=\pi_1(T^3\smallsetminus  G^*)$;
the assumption on $H_2$ is satisfied trivially since $H_2(T^3\smallsetminus G^*)=0$. 
It follows that the map $J\longrightarrow H$ (where $H$ now denotes the first homology of the preimage of $ G^*$ in $\R^3$) is injective, since $J=F^{(1)}/F^{(2)}$, $H=B^{(1)}/B^{(2)}$. 

This approach does not work for links in place of the $1$-complex $ G^*$ because $H_2(T^3\smallsetminus  L)$ is non-zero, and for a given link $L$ it is a non-trivial problem to determine what the second homology is of the image on $\pi_1$ of an arbitrary map $G\longrightarrow T^3\smallsetminus L$.

\section{Powers of the augmentation ideal and the lower central series} \label{sec: exact seq}

The result of this section, in conjunction with Lemmas \ref{augmentation ideal lemma},  \ref{links lemma} and their analogues for $T^3$ discussed in section \ref{3-torus section}, 
completes the proof of theorem \ref{k-filling thm}.
Let $i'\co G'\longrightarrow T^3\smallsetminus L$ denote any spine homotopic to the standard spine $G$, where $L$ is a link whose components are all essential in $\pi_1 (T^3)$, as in section \ref{3-torus section}. The analysis below also applies to the relative case $T^2\times I$, considered in section \ref{FM and powers}. 

We start by recalling the notation  and summarizing basic consequences of the topological setup.  $F$ denotes $\pi_1 (G')$, the free group on three generators, and
$$K:= {\rm ker}\, [\, \pi_1(T^3\smallsetminus L)\longrightarrow \pi_1(T^3)\, ] $$
is isomorphic to $\pi_1({\mathbb R}^3\smallsetminus \widetilde L)$, a free group.
$J'$ denotes the first homology of the preimage $\widetilde G'$ of $G'$ in $\R^3$ and  $H$ denotes $H_1({\mathbb R}^3\smallsetminus \widetilde L)$,
$$
J'\cong F_2/[F_2, F_2], \; \, H\cong K/[K,K].$$

$J'$ and $H$ are considered as modules over $\Z[\Z^3]$, and $I^k$ denotes the $k$-th power of the augmentation ideal $I$ of $\Z[\Z^3]$.
The following statement is the main result of this section, relating the filtrations of $J'$, $H$ in terms of powers of the augmentation ideal and  the lower central series of $F$, $\pi_1(T^3\smallsetminus L)$.

\begin{lemma} \label{augm lcs} {\sl Fix $k\geq 2$, and 
suppose $I^j J'/I^{j+1} J'\longrightarrow  I^j H/I^{j+1}  H$ is injective for all $0 \leq j\leq k-2$.
Then the map
$$F/F_{k+1} \longrightarrow  \pi_1(T^3\smallsetminus L)/\pi_1(T^3\smallsetminus L)_{k+1},$$
induced by the inclusion $i'\co G'\longrightarrow T^3\smallsetminus L$,
 is injective. 
}
\end{lemma}

It is convenient to introduce the following notation, so that the proof involves isomorphisms rather than injections. Consider
\begin{equation} \label{pi eq} 
\pi:={\rm image}\; [ \,  F \overset{i'_*}\longrightarrow \pi_1(T^3\smallsetminus L)\,].
\end{equation}
 
Since the composition $F\longrightarrow \pi\longrightarrow \pi_1 (T^3)$ is an isomorphism on $H_1$, $H_1 F\longrightarrow H_1\pi$ is an isomorphism. 
Note that for $k>1$, $\pi_k$ is contained in $K$, and so is a free group.
Consider  $$\overline H\,: =\, {\rm image}[\, J'\overset{i'_*}\longrightarrow H\, ]\; \cong\; (\pi\, \cap \, K) /[\pi\, \cap \, K,\pi\, \cap \, K]. $$

Observe that $\pi\, \cap\, K=\pi^{}_2$. Indeed, as noted above $\pi_2\subset K$. To prove the opposite inclusion, consider $g\in \pi\cap K$. Then $g=i'_*(f)$ for some $f\in F$, see diagram (\ref{commutator subgroup}). Since $f\in {\rm ker}\, [\, F\longrightarrow \pi_1 (T^3)\, ]\, =F^{}_2$, it follows that $g\in \pi^{}_2$.
\begin{equation} \label{commutator subgroup}
\centering
\begin{tikzcd}
1 \arrow{r} & \pi\, \cap \, K\arrow{r} & \pi \arrow{r} & \pi_1 (T^3)\arrow{r} & 1 \\
& F \arrow[ur,twoheadrightarrow,"i'_*"'] & & & 
\end{tikzcd}
\end{equation}

Therefore $\;\overline H$ could also be defined as $\overline H=\pi^{}_2/[\pi^{}_2,\pi^{}_2]$.

{\em Proof of lemma \ref{augm lcs}.}
The proof is by induction on $k$. We add the base case $k=1$ where there is no assumption on the link, and $H_1 F\longrightarrow H_1\pi$ is an isomorphism (the conclusion of the lemma when $k$ is set to $1$) as discussed above. Suppose the statement of the lemma holds for $k-1$, so the inductive assumption is that $F/F_k\longrightarrow \pi/\pi_k$ is an isomorphism. The overall strategy is motivated by the proof of the Stallings theorem \cite{S65}; in particular we use the isomorphism $H_2(F/F_k)\cong H_2(\pi/\pi_k)$ which is a consequence of the inductive assumption. It is worth noting again that both $F^{}_2$ and ${\pi}^{}_2$ are free groups but the map $F^{}_2\longrightarrow {\pi}^{}_2$ is not an isomorphism on $H_1$. Being an isomorphism on $H_1$ is equivalent to $J'/I^k J'\cong  \overline H/I^k \overline H$ for {\em all} $k$.
Rather the lemma has a weaker assumption, $J'/I^k J'\cong  \overline H/I^k \overline H$ for {\em some} fixed $k$. 

Therefore the main technical ingredient in the proof, the diagram of exact sequences (\ref{cd}), is formulated with the goal of comparing the lower central series to powers of the augmentation ideal and it differs from the one used in the Stallings theorem in \cite{S65}. 

We start by setting up the relevant short exact sequences.
Let $\phi_k$ denote the inclusion $F_k\subset F^{}_2$ composed with the quotient map $F^{}_2 \longrightarrow F^{}_2/[F^{}_2,F^{}_2]$, and consider its kernel:
\begin{equation} \label{SES1}
1\longrightarrow F_k\cap [F^{}_2,F^{}_2] \longrightarrow F_k\overset{ \phi_k}\longrightarrow F^{}_2/[F^{}_2,F^{}_2]
\end{equation} 
 
Denote the generators of $F$ by $x,y,z$; the same letters will denote the covering translations of ${\mathbb R}^3$. To relate the map $\phi_k$ to the geometric discussion of the ``jungle gym'' $J$ in the previous section, consider a basic example, the triple commutator
$[[x,y],z]\in F^{}_3$. The map $\phi_k$ is implemented by first expanding 
$[[x,y],z]=[x,y]\cdot([x,y]^{-1})^z$. The first factor is mapped to the boundary of the plaquette $P_z$, figure \ref{fig:lines}. The second factor is mapped to the boundary of this plaquette with the opposite orientation and shifted one unit up, $-zP_z$. So $\phi_3([[x,y],z])=(1-z)P_z$.
An arbitrary element of $F^{}_3$ expands as a product of conjugates of elements of the form $[g_1,g_2]^{g_3}\cdot([g_1,g_2]^{-1})^{g_4}$, where each $g_i\in F$. The map $\phi_3$ takes $[g_1,g_2]$ to a cycle $c$ in $\widetilde G'$, and 
\begin{equation} \label{commutators eq}{\phi}_3([g_1,g_2]^{g_3}\cdot([g_1,g_2]^{-1})^{g_4}\; = \; (g_3-g_4)\cdot c\, \in \, IJ'.
\end{equation}

It follows that the map $\phi_3$ surjects onto $IJ'$.
The analogous statements hold for $F_k$: the  map $\phi_k$ takes a basic commutator $[\ldots[[x_1,x_2],x_3],\ldots,x_k]$, where $x_1,\ldots, x_k\in\{x,y,z\}$, $x_1\neq x_2$,
to the plaquette determined by $[x_1,x_2]$, multiplied by 
 $(1-x_3)\ldots(1-x_k)\in I^{k-2}$. More generally,
the image of $\phi_k$ is in $I^{k-2}J'\subset J'$, and moreover any element of $I^{k-2}J'$ is in the image of $\phi_k$. Consider the exact sequence (\ref{SES1}) for $k, k+1$:

 \begin{equation} \label{cd0} \begin{tikzcd} \large
 1\arrow{r} & F_k\cap [F_2,F_2] \arrow{r}  & F_k\arrow[r, "\phi_k"]  & I^{k-2}J' \arrow{r}  & 1 \\
1\arrow{r} & F_{k+1}\cap [F_2,F_2] \arrow{r} \arrow[hookrightarrow]{u}  & F_{k+1}\arrow[r, "\phi_{k+1}"] \arrow[hookrightarrow]{u} & I^{k-1} J' \arrow{r} \arrow[hookrightarrow]{u} & 1
\end{tikzcd}
\end{equation}

The quotients of the respective groups form the short exact sequence

\begin{equation} \label{SES4}
 1\longrightarrow \frac{F_k\cap [F_2,F_2]}{F_{k+1}\cap [F_2,F_2]} \longrightarrow \frac{F_k}{F_{k+1}}\longrightarrow \frac{I^{k-2}J'}{I^{k-1}J'} \longrightarrow 1
\end{equation} 

Recall that for $k>1$, $\pi_k$ is contained in $K$, and is a free group.  
We would like to compare the short exact sequence (\ref{SES4}) for $F$ to the corresponding one for $\pi$.
The sequence for $\pi$ is obtained by starting with the map $\pi_k \longrightarrow \pi^{}_2/[\pi^{}_2,\pi^{}_2]$ as in (\ref{SES1}). The same analysis as above, for example equation (\ref{commutators eq}) with $J'$ replaced with $\overline H$, shows that $\pi_k$ surjects onto $I^{k-2} \overline H$. 
The vertical maps in the following diagram, relating the short exact sequences (\ref{SES4}) for $F, \pi$, are induced by the inclusion $G'\subset T^3\smallsetminus L$.
\tikzcdset{every label/.append style = {scale=1,yshift=0.6ex},
every matrix/.append style={nodes={font=\large}}
} 
\begin{equation} \label{cd} \begin{tikzcd} \large
 1\arrow{r} & \frac{F_k\, \cap \, [F^{}_2,F^{}_2]}{F_{k+1}\, \cap\,  [F^{}_2,F^{}_2]} \arrow{r} \arrow[d,"\alpha"] & \frac{F_k}{F_{k+1}}\arrow{r} \arrow[d, "\beta"] & \frac{I^{k-2}J'}{I^{k-1}J'} \arrow{r} \arrow[d, "\gamma"] & 1 \\
1\arrow{r} & \frac{\pi_k\, \cap \, [\pi^{}_2,\pi^{}_2]}{\pi^{}_{k+1}\, \cap\,  [\pi^{}_2,\pi^{}_2]} \arrow{r}  & \frac{\pi_k}{\pi^{}_{k+1}}\arrow{r} & \frac{I^{k-2}\overline H}{I^{k-1}\overline H} \arrow{r} & 1
\end{tikzcd}
\end{equation}
The inductive assumption is that $F/F_k\longrightarrow \pi/\pi_k$ is an isomorphism, and the goal is to show that the middle vertical map $\beta$ is an isomorphism, to propagate the inductive step. The assumption in the statement of the lemma, based on the link $L$, implies that the map $\gamma$ is an isomorphism. 
The goal is to show that ${\alpha}$ is an isomorphism; then it will follow that $\beta$ is an isomorphism too, concluding the inductive step. 

Recall that $F\longrightarrow \pi$ is surjective and $F/F_k\longrightarrow \pi/\pi_k$ is assumed to be an isomorphism. It follows that $F^{}_2/F_k\longrightarrow {\pi}^{}_2/\pi_k$ is an isomorphism, and so is $H^{}_2(F^{}_2/F_k)\longrightarrow H^{}_2({\pi}^{}_2/\pi_k)$. Consider Hopf's characterization of $H^{}_2$ of a group $ {\rm Free}/{\rm R}$:
{\small
\begin{equation} \label{H2}
H^{}_2({\rm F}/{\rm R})
%H_2\left(\frac{\rm Fr}{\rm R}\right)
\; \cong \;  \frac{ {\rm R}\cap [{\rm F}, {\rm F}]}{[{\rm F}, {\rm R}]}
\end{equation}
}
Apply this to $H^{}_2(F^{}_2/F_{k})\cong H^{}_2({\pi}^{}_2/{\pi}_k)$ to get
\begin{equation} \label{H2comp}
\frac{F_k\, \cap \, [F^{}_2,F^{}_2]}{  [F_k,F^{}_2]}\; \cong \; 
\frac{\pi_k\, \cap \, [\pi^{}_2,\pi^{}_2]}{  [\pi_k,\pi^{}_2]}
\end{equation}
In the Stallings' proof of his theorem, the second homology groups are terms in the 5-term exact sequence \cite[Theorem 2.1]{S65}; in our context
the relation is not as immediate: the map $\alpha$ in (\ref{cd}) and the isomorphism (\ref{H2comp}) are related by the short exact sequences

 \begin{equation} \label{cd1} \begin{tikzcd} \large
 1\arrow{r} & 
 \frac{F_{k+1}\, \cap \, [F_2,F_2]}{  [F_k,F_2]} \arrow{r} \arrow{d} &
 \frac{F_k\, \cap \, [F_2,F_2]}{  [F_k,F_2]} \arrow{r} \arrow[d, "\cong"] &
 \frac{F_k\, \cap \, [F_2,F_2]}{F_{k+1}\, \cap\,  [F_2,F_2]} \arrow{r} \arrow[d, "\alpha"] & 1 \\
1\arrow{r} & 
\frac{\pi^{}_{k+1}\, \cap \, [\pi_2,\pi_2]}{  [\pi_k,\pi_2]} \arrow{r}  &
\frac{\pi_k\, \cap \, [\pi_2,\pi_2]}{  [\pi_k,\pi_2]}\arrow{r} &
\frac{\pi_k\, \cap \, [\pi_2,\pi_2]}{\pi^{}_{k+1}\, \cap\,  [\pi^{}_2,\pi^{}_2]} \arrow{r} & 1
\end{tikzcd}
\end{equation}
 
The middle vertical map is an isomorphism by (\ref{H2comp}). It follows that the left vertical map is injective. Since $F\longrightarrow \pi$ is surjective, the left map is in fact an isomorphism. It follows that $\alpha$ is an isomorphism.
\qed
 
We are in a position to show how the statements established in Sections \ref{FM and powers} - \ref{sec: exact seq} imply the proof of Theorem \ref{k-filling thm}.

{\em Proof of Theorem \ref{k-filling thm}.}
As discussed in the paragraph following equation \eqref{pi eq}, any link (for example, the empty link) is $2$-filling.
Consider $k\geq 3$.
According to Lemma \ref{links lemma}, there exists a link $L_{k-3}\subset T^2\times I$ such that the map induced by the inclusion of the standard spine $G$ into the link complement, $ I^j J/I^{j+1} J\longrightarrow  I^j H/I^{j+1}  H$, is injective for all $1 \leq j\leq k-3$. Its analogue for $T^3$ is established in section \ref{3-torus section}. Lemma \ref{augmentation ideal lemma} shows that then for any spine $G'$ in the link complement, $ I^j J'/I^{j+1} J'\longrightarrow  I^j H/I^{j+1}  H$ is injective for all $1 \leq j\leq k-3$. Let $M$ denote either $T^2\times I$ or $T^3$. Now it follows from Lemma \ref{augm lcs} that the map of fundamental groups modulo the $k$-th term of the lower central series $F/F_{k} \longrightarrow  \pi_1(M\smallsetminus L)/\pi_1(M\smallsetminus L)_{k}$, induced by the inclusion $ G'\longrightarrow M\smallsetminus L$, is injective, thus completing the proof of the theorem.
\qed

\subsection{Discussion and questions}\label{sec: question}

For each $k>1$, the combination of lemmas \ref{augmentation ideal lemma}, \ref{links lemma} and \ref{augm lcs} gives a $k$-filling link $L_k\subset T^3$: for any $1$-spine $G'$ of $T^3$ which is disjoint from $L_k$, $$F/F_k \longrightarrow   \pi_1(T^3\smallsetminus L_k)/\pi_1(T^3\smallsetminus L_k)_k$$ 

is injective. (As above, $F$ denotes $\pi_1(G')$.) Denoting by $\pi$ the image of $F$ as in (\ref{pi eq}), this means $F/F_k\cong\pi/\pi_k$.
It is interesting to note that by Dwyer's theorem  (see section \ref{background sec}) this implies that $H_2(\pi)$ is contained in the $(k-1)$st term of the Dwyer filtration $\phi_{k-1}(\pi)$. This is true for {\em any} map $G'\subset T^3\smallsetminus L_k$, a fact that seems quite non-trivial to prove directly since $H_2(T^3\smallsetminus L_k) \nsubseteq {\phi}_3(T^3\smallsetminus L_k)$.

The cardinality of the links $L_k$ grows linearly with $k$. 
In this paper we considered a special collection of links $L$ in $T^3$, ensuring that the complement in the universal cover has free fundamental group.
Using this approach in the general case would require a non-trivial extension of the group-theoretic analysis developed above. 

To illustrate the last point, we present a specific example, a ``dense chain mail link'', for which the property of being filling is not known to us.  
Consider the unit cubical $T^3$, and take an $\epsilon$-net of points in $T^3\times S^2$, the Grassmanian of all tangent $2$-planes to $T^3$, for some small $\epsilon$.  Now define the link L to be the link of circles of radius $1/10$ in $T^3$ centered at and oriented by each point of the net. By perturbing the net we may assume L is an embedded link. Now suppose $G$ is a spine of $T^3$ embedded disjointly from $L$. 
Each component of $L$ bounds a disk in $T^3$; of course the disks for nearby components intersect. Thin finger move homotopies can be used to push $G$ off these disks, ensuring trivial linking number with $L$. (Generically these finger moves will acquire intersections with other disks, but the new intersections come in $\pm 1$ pairs and thus do not contribute to linking numbers.) The linking number gives information about the induced map on homology but the question of injectivity modulo the third and higher terms of the lower central series is open.

{\bf Acknowledgements.} We would like to thank Danny Calegari for discussions on the still open question of the existence of filling links for closed $3$-manifolds of higher rank.

VK was supported in part by the Miller Institute for Basic Research in Science at UC Berkeley, Simons Foundation fellowship 608604, and NSF Grant DMS-1612159.

\addresseshere

\newpage
\appendix

\setcounter{footnote}{0}

\section{Filling links in $3$-manifolds of rank $2$}
\centerline{\bf by Christopher J. Leininger\footnote{Supported in part by NSF grant 
DMS-2106419}  and Alan W. Reid\footnote{Supported in part by NSF grant DMS-1812397} }

In this appendix we prove the following result which provides the first examples of closed orientable $3$-manifolds with Heegaard genus $>1$ that contain a filling link. To state this, recall that the {\em rank} of a finitely generated
group $\Gamma$ is the minimal cardinality of a generating set for $\Gamma$. When $\Gamma=\pi_1(M)$ and $M$ is a compact $3$-manifold, we define the {\em rank of $M$} to be the rank of $\Gamma$. 

\begin{theorem}
\label{rank2}
Let $M$ be a closed orientable $3$-manifold of rank $2$. Then $M$ contains a filling hyperbolic link.\end{theorem}

For convenience, we will discuss what lies behind the crucial feature that we exploit to exhibit filling links using the assumption of rank $2$; namely a classical result from $3$-manifold topology (see \cite[Theorem VI.4.1]{JS}) which affords a classification of $2$-generator subgroups of the fundamental group of a compact atoroidal, irreducible $3$-manifold (see also \S \ref{proofr2} below). Very briefly, let $X$ be such a 3-manifold, $H\subset \pi_1(X)$ a $2$-generator subgroup, and $Y_H\rightarrow X$ the cover corresponding to $H$ with compact core $C_H$. 
As an illustration of the main part of the argument, assume that $\pi_1(X)$ is freely indecomposable and $\partial C_H$ is non-empty and contains no $2$-sphere components. Then the $2$-generator assumption, together with the fact that the
first Betti number of a $3$-manifold is at least half the first Betti number of its boundary, implies that $\partial C_H$ (which may be disconnected) has genus $1$ or $2$. In the case when $\partial C_H$ is connected of genus $2$, 
it follows that $b_1(C_H)=2$, and a standard argument in $3$-manifold topology now shows:
$$- 1 = \chi(\partial C_H)/2 = \chi(C_H) = b_2(C_H) - b_1(C_H) + b_0(C_H) =  b_2(C_H) - 1.$$
It follows that $b_2(C_H)=0$, and one can now deduce in this case that $H$ is free of rank $2$, contradicting the assumption that $H$ is indecomposable. The remainder of the argument is completed by
analyzing the case when $\partial C_H$ consists of tori; either one or two incompressible tori, or when $\partial C_H$ contains a compressible torus. The upshot is a limited set of possibilities for what $C_H$, and hence $H$ can be, and these are listed in \S \ref{proofr2}. 

Before commencing with the proof of Theorem \ref{rank2} we list some corollaries.
Closed orientable $3$-manifolds of Heegaard genus $2$ are examples of manifolds covered by Theorem \ref{rank2}, so an immediate corollary is:

\begin{corollary}
\label{genus2}
Let $M$ be a closed orientable $3$-manifold such that $M$ has Heegaard genus $2$. Then $M$ contains a filling hyperbolic link.\end{corollary}

Note that there there are closed orientable $3$-manifolds $M$ for which $\pi_1(M)$ has rank $2$, but the Heegaard genus is $3$
(see \cite{BZ}). However, an interesting case of Corollary \ref{genus2} is the following. 

Using Perelman's resolution of the Geometrization Conjecture \cite{KL}) it is known that all closed orientable $3$-manifolds with non-trivial finite fundamental  have Heegaard genus $1$ or $2$. To see this, the resolution of the Geometrization Conjecture proves that a closed $3$-manifold $M$ with finite fundamental is covered by $S^3$, in which case $M$ is a Seifert fibered space. Assuming that $M$ is not $S^3$, or a Lens Space, then $M$ is a Seifert fibered space
over $S^2$ with three exceptional fibers (see \cite[VI.11]{Ja}), from which a genus $2$ Heegaard splitting may be constructed directly (see \cite[Proposition 1.3]{BZ} for example).
The existence of filling links in manifolds of genus $1$ was noted after the proof of Proposition 2,  so Theorem \ref{rank2} now shows.

\begin{corollary}
\label{finite}
Let $M$ be a closed orientable $3$-manifold such that $\pi_1(M)$ is finite and non-trivial. Then $M$ contains a filling link.\end{corollary}

\begin{remark} An alternative proof of Corollary \ref{finite} bypassing the use of Heegaard genus is the following. Using \cite[Section 3]{Mil}, one obtains a classification of finite groups that act freely on $S^3$. Perelman's 
resolution of the Geometrization Conjecture eliminates the one class of finite groups from \cite{Mil} that are not subgroups of $\SO(4)$.  The five families of non-cyclic subgroups of $\SO(4)$ can all be seen to have rank $2$.\end{remark}

\subsection{Proof of Theorem \ref{rank2}.} 
\label{proofr2}
We begin with a preliminary remark. The hypothesis that $M$ is closed, orientable and $\pi_1(M)$ has rank $2$ implies that $\pi_1(M)$ is non-Abelian. The reason is this:  From \cite[Theorem 9.13]{He}
the only abelian groups occurring as the fundamental group of a closed orientable 3-manifold are $\Z$, $\Z/n\Z$, and $\Z^3$, and these are all excluded by the rank hypothesis.

Let $L\subset M$ be a hyperbolic link with at least $3$ components.  To find such a link, one may start with any link $L'$ in $M$ with at least $2$ components, then appealing to \cite[Corollary 6.3]{M82} we may choose a knot $K$ in $M \setminus L'$ with hyperbolic complement, and set $L = L' \cup K$.

As in the introduction $G$ will be a spine of $M$, with $G\cap L=\emptyset$.  By definition, $\pi_1(G)$ surjects onto $\pi_1(M)$. We need to prove that the induced homomorphism $\pi_1(G)\rightarrow \pi_1(M\setminus L)$ is injective.  Let $H$ be
be the image group of this homomorphism. Recall that $\pi_1(G)$ is free of rank $2$, and so 
the result will follow from the Hopfian property for free groups once we establish that $H$ is free of rank $2$.  We remind the reader that the Hopfian property for a free group $F$ asserts that every epimorphism $F\rightarrow F$ is an isomorphism.

We argue as follows. Appealing to the classical $3$--manifold result sketched in the previous section, namely \cite[Theorem VI.4.1]{JS}, observe that since $M\setminus L$ is hyperbolic the possibilities for $H$ are:

\medskip

\begin{enumerate}
\item $H$ is free of rank $2$, or
\item $H$ is free abelian of rank $\leq 2$, or
\item $H$ finite index in $\pi_1(M\setminus L)$.
\end{enumerate}

\medskip

If $H$ is free abelian of rank $1$ or $2$, then the image in $\pi_1(M)$ via the homomorphism induced by inclusion $M\setminus L \to M$ is the quotient of a free abelian group of rank at most $2$, and is in particular abelian.  
Since $G$ is a spine, $H$ surjects onto $\pi_1(M)$, 
which contradicts the fact noted above that $\pi_1(M)$ is non-abelian.  Therefore, case (2) is impossible.
Case (3) can also be eliminated as follows.  As noted above, the first Betti number of a $3$-manifold is at least half the first Betti number of its boundary, and so we deduce from this that $\H^3/H$ has at most two cusps.  On the other hand,
by construction, $M\setminus L$ has at least $3$ cusps, and $\H^3/H\rightarrow M\setminus L$ is a finite sheeted cover, a contradiction. We therefore conclude that the only possibility for $H$ is that it is free of rank $2$ and the proof is complete.\qed

\bigskip

{\small{
\sc{ Department of Mathematics\\
Rice University\\
Houston, TX 77005, USA}

{\em E-mail address}: {\tt c.j.leininger95@rice.edu}\\
{\em E-mail address}: {\tt alan.reid@rice.edu}
}}


\begin{thebibliography}{F82b}




\bibitem[Bi58]{Bing} R.H. Bing, {\em Necessary and sufficient conditions that a 3-manifold be $S^3$}, Ann. of Math. 68 (1958), 17-37.

\bibitem[BZ84]{BZ} M.~Boileau and H.~Zieschang, {\em Heegaard genus of closed orientable $3$-manifolds}, Invent. Math. {\bf 76} (1984), 455-468.


\bibitem[CG83]{CG83}  A.J. Casson and C.McA. Gordon, {\em A loop theorem for duality spaces and fibred ribbon knots} Invent. Math. 74 (1983), 119-137.

\bibitem[CH08]{CH08} T.D. Cochran and S. Harvey, {\em Homology and derived series of groups. II. Dwyer's theorem}, Geom. Topol. 12 (2008), 199-232.

\bibitem[CO98]{CO98} T.D. Cochran and K. Orr, {\em Stability of lower central series of compact 3-manifold groups}, Topology 37 (1998), 497-526.


\bibitem[DK92]{DK92} G. Duchamp, and D. Krob, {\em The lower central series of the free partially commutative group},
Semigroup Forum 45 (1992), 385-394.

\bibitem[Dw75]{D75} W.G. Dwyer,
{\em Homology, Massey products and maps between groups},
J. Pure Appl. Algebra 6 (1975), 177-190.

\bibitem[Fr20]{F20}  M. Freedman, {\em Controlled Mather-Thurston theorems}, arXiv:2006.00374

\bibitem[FK20]{FK20} M. Freedman and V. Krushkal, {\em Engel groups and universal surgery models}, J. Topol. 13 (2020), 1302-1316.

\bibitem[FQ90]{FQ90} M. Freedman and F. Quinn, Topology of $4$-manifolds, Princeton Mathematical Series, 39. Princeton University Press, Princeton, NJ, 1990.

\bibitem[He76]{He} J.~Hempel, {\em 3-Manifolds}, Annals of Math. Studies {\bf 86}, P.U.P. (1976).



\bibitem[Ja80]{Ja} W.~Jaco, {\em Lectures on Three-Manifold Topology}, CBMS Regional Conference Series in Mathematics, {\bf 43} A. M. S. Providence, R.I., (1980).

\bibitem[JS79]{JS}  W.~Jaco and P.~B.~Shalen, {\em Seifert Fibered Spaces in $3$-Manifolds}, Mem. Amer. Math. Soc. {\bf 21} , no. 220 (1979), 192pp.



\bibitem[KL08]{KL} B.~Kleiner and J.~Lott, {\em Notes on Perelman's papers}, Geometry and Topology {\bf 12} (2008), 2587-2855.


\bibitem[MKS66]{MKS66} W. Magnus, A. Karrass, and D. Solitar, Combinatorial group theory: Presentations of groups in terms of generators and relations. Interscience Publishers [John Wiley \& Sons, Inc.], New York-London-Sydney 1966.

\bibitem[Me18]{M18} G. Meigniez, {\em Quasi-complementary foliations and the Mather-Thurston theorem}, arXiv:1808.02377

\bibitem[Mi57]{Mil} J.~Milnor, {\em Groups which act on $S^n$ without fixed points}, Amer. J. Math. {\bf 79} (1957), 623-630.


\bibitem[My82]{M82} R. Myers, 
{\em Simple knots in compact, orientable 3-manifolds},
Trans. Amer. Math. Soc. 273 (1982), 75-91.

\bibitem[St21]{Stagner} W. Stagner, {\em $3$-manifolds of rank $3$ have filling links}, arXiv:2110.02936


\bibitem[St65]{S65} J. Stallings, {\em Homology and central series of groups}, J. Algebra 2 (1965), 170-181.

\bibitem[St75]{S75}  J. Stallings, 
{\em Quotients of the powers of the augmentation ideal in a group ring}, Knots, groups, and 3-manifolds,  Ann. of Math. Studies, No. 84, Princeton Univ. Press, Princeton, N.J., 1975, pp. 101–118.

\end{thebibliography}
\end{document}